# A Fourier series solution for the longitudinal vibrations of a bar with viscous boundary conditions at each end

Vojin Jovanovic
Systems, Integration & Implementation
Smith Bits, A Schlumberger Co.
1310 Rankin Road
Houston, TX 77073
e-mail: fractal97@gmail.com

## Abstract

This paper presents the generalized Fourier series solution for the longitudinal vibrations of a bar subjected to viscous boundary conditions at each end. The model of the system produces a non-self-adjoint eigenvalue problem which does not yield a self-orthogonal set of eigenfunctions with respect to the usual inner product. Therefore, these functions cannot be used to calculate the coefficients of expansion in the Fourier series. Furthermore, the eigenfunctions and eigenvalues are complex-valued. The eigenfunctions can be utilized if the space of the wave operator is extended and a suitable inner product is defined. It is further demonstrated that the series solution contains the solutions for free-free, fixed-damper, fixed-fixed and fixed-free bar cases. The presented procedure is applicable in general to other problems of this type. As an illustration of the theoretical discussion, the results from numerical simulations are presented.

## 1. Introduction

This paper presents an investigation of the problem of longitudinal vibrations of a bar with viscous boundaries at each end. Motivation and purpose for this study is due to the intricacies that arise when viscous dampers are used to control the displacement of large structures such as buildings and bridges. Such control is especially needed in cases when mitigating hazardous effects of earthquakes is desired. During an earthquake, a finite quantity of energy is introduced into a structure in which the level of damping is typically very low. If there were no damping, vibrations would continue indefinitely. Since there is always some level of damping the vibrations cease. However, the response of the structure to input energy can be improved if energy is absorbed not by the structure itself but by a viscous boundary. Unfortunately, such a boundary introduces the boundary conditions that drastically change the mathematical nature of the problem. The eigenvalues of the underlying mathematical problem become complex-valued and the eigenmodes become non-orthogonal with respect to the usual inner product. Such a situation requires more general procedures which are presented in this paper.

While from a computational point of view the methodology presented here may not be easier than applying the finite element method (FEM), it is crucial to properly interpret FEM results. Since the eigenvalues are complex-valued, which usually does not occur in practice, a standard application of FEM using commercial software could lead to wrong interpretation of the stability of the system. This is due to spurious eigenvalues in FEM and is discussed in this paper.

Historically, attempts to treat non-self-adjoint boundary conditions in the engineering literature are sparse perhaps because the vast majority of cases lead to self-adjoint mathematical description. Nevertheless, several references exist such as [1], [2], [3] and[4] which provide certain insights into non-self-adjoint nature of longitudinal vibrations of a bar with one viscous boundary. Reference [5] discusses a non-self-adjoint problem of acoustic ducts which is equivalent to the vibration of free-fixed bar. Of these references [3] appears to be the first one to provide the complete series solution for that problem. A previous attempt reported in [4] utilized a non-standard treatment for a decoupling the equation of motion and provided a response only to a harmonic driving force. However, imposing an additional viscous boundary condition at the other end of the bar, introduces a rigid body mode to the system. This in turn introduces non-intuitive behavior as discussed in [6]. Reference [5] uses the state space approach in a way different from the one in this paper. The authors cast the problem into the first order form and proceed to show that a new inner product orthogonalizes the eigenmodes and implicitly decouples the system. Dealing with this new inner product is somewhat cumbersome and details are given on going back from the state space formulation to calculating the displacement function. There is also no discussion of the rigid body mode. We demonstrate, however, that the usual inner product can be retained if one makes use of the eigenmodes of the dual problem. Biorthogonality makes the determination of the expansion coefficients straightforward with all integrations carried over the original domain (without the artificial mirroring in reference [5]), despite the fact that we have a more general boundary condition on the left side of the bar. In this way the state space system is decoupled explicitly, and we supply all the details for computing the expansion coefficients and the complete response of the bar. We also discuss the rigid body mode and the limiting cases of our boundary conditions.

We organize the presentation as follows. We first state the problem and provide its mathematical description. From it we derive a boundary value problem (BVP) and determine eigenfunctions and eigenvalues which turn out to be complex-valued. We then analyze the differential operator and recast the BVP into a state-space self-adjoint form which yields bi-orthogonal eigenvectors. We then perform a usual modal analysis step together with the Laplace transform and write a response of the system in terms of eigenfunction expansion. We proceed further with qualitative analysis of the solution to verify the known eigenfunctions for free-free, fixed-damper, fixed-fixed and fixed-free bar cases. Finally, we provide numerical illustrations and discuss the efficiency of the finite element approximation relative to the analytical solution. We also point out to, at this time, an unexplained phenomenon pertaining to FEM's distortion of the stability regions making a physically stable system appear unstable.

## 2. Problem statement

We begin with the problem statement of the system. Figure 1 depicts a bar suspended by a damper at each end but free to move horizontally. The symbols $\rho$, $A_0$ and $E$ represent the density of the bar, constant cross-sectional area and the modulus of elasticity respectively. The coefficients $c_1$ and $c_2$ represent damping at each end of the bar.

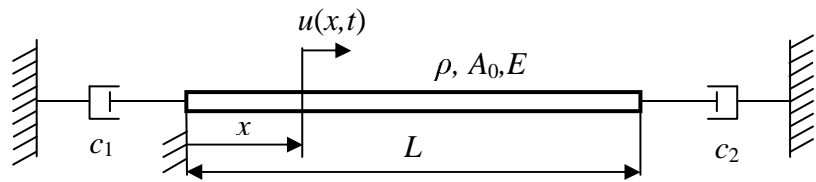

**Figure 1. A bar with two viscous boundaries.**

A wave traveling along the bar must satisfy the following equation of motion - the wave equation -

$$u_{tt}(x,t) = c^2 u_{xx}(x,t) + p(x,t), \quad 0 < x < L, \tag{1}$$

with its associated boundary

$$u_x(0,t) = \frac{h_1}{c} u_t(0,t) \quad \text{and} \quad u_x(L,t) = -\frac{h_2}{c} u_t(L,t) \tag{2}$$

and initial conditions

$$u(x,0) = f(x) \quad \text{and} \quad u_t(x,0) = g(x). \tag{3}$$

Subscripts $x$ and $t$ denote partial derivates with respect to space and time respectively, while $c^2 = \frac{E}{\rho}$, $h_1 = \frac{c_1 c}{E A_0}$ and $h_2 = \frac{c_2 c}{E A_0}$ where $c$ is the speed of a wave along the bar and $h_1$ and $h_2$ are non-dimensional coefficients. Functions $f(x)$ and $g(x)$ are assumed to be real square-integrable functions representing initial displacement and velocity of the bar. They are described relative to the fixed coordinate system. $p(x,t)$ is the force per unit mass distributed along the length of the bar while $u(x,t)$ is the displacement coordinate of a point of the bar at the initial location $x$.

## 3. Eigenvalues and eigenfunctions

The wave equation can be solved using the separation of variables method by assuming

$$u(x,t) = \varphi(x) e^{\lambda t}. \tag{4}$$

Substituting this into the homogeneous wave equation yields the following BVP

$$\varphi''(x) - \frac{\lambda^2}{c^2} \varphi(x) = 0 \tag{5}$$

$$\varphi'(0) = \frac{h_1}{c} \lambda \varphi(0)$$

$$\varphi'(L) = -\frac{h_2}{c} \lambda \varphi(L).$$

To solve this eigenvalue problem we set

$$\varphi(x) = A \sinh\left(\frac{\lambda}{c} x\right) + B \cosh\left(\frac{\lambda}{c} x\right) \tag{6}$$

where $A$ and $B$ are arbitrary (possibly complex) constants and where hyperbolic functions are utilized mostly for convenience. We note that there is no any assumption on the spectral parameter $\lambda$ which can be either real or complex. The boundary condition at $x = 0$ reveals that either $\lambda = 0$ or $A = h_1 B$. The former imposes the rigid body eigenvalue $\lambda_r = 0$ which corresponds to its associated eigenfunction $\varphi_r = A_r = const$. The latter imposes an eigenfunction

$$\varphi(x) = B \left( \cosh\left(\frac{\lambda}{c} x\right) + h_1 \sinh\left(\frac{\lambda}{c} x\right) \right) \tag{7}$$

for which the associated eigenvalues (all different from zero) can be obtained from the boundary condition at $x = L$ as

$$\tanh\left(\frac{\lambda L}{c}\right) = -\frac{h_1 + h_2}{1 + h_1 h_2}. \tag{8}$$

Eq. (8) can be recast into the form

$$e^{2\frac{\lambda L}{c}} = \frac{(1-h_1)(1-h_2)}{(1+h_1)(1+h_2)}. \tag{9}$$

The solution to Eq. (9) is

$$\lambda_n = \frac{c}{2L}\left[\ln\left|\frac{(1-h_1)(1-h_2)}{(1+h_1)(1+h_2)}\right| + i\left(Arg\left(\frac{(1-h_1)(1-h_2)}{(1+h_1)(1+h_2)}\right) + 2n\pi\right)\right], \quad n \in \mathbb{Z} \tag{10}$$

Therefore, the eigenvalues for BVP are $\lambda_r, \lambda_0, \lambda_{\pm 1}, \lambda_{\pm 2},....$ while the eigenfunctions are $\varphi_r, \varphi_0, \varphi_{\pm 1}, \varphi_{\pm 2},.....$ The eigenvalues and eigenfunctions for special cases $h_1 = \pm 1$, $h_2 = \pm 1$, $h_2 = -h_1$ and $h_2 = -1/h_1$ are discussed in [6].

## 4. Adjoint operator and Fourier's expansion

It can be easily verified that the eigenfunctions Eq.(7) which satisfy the BVP are not mutually orthogonal with respect to the usual standard inner product. This property is required for obtaining the eigenfunction expansion of the exact solution for the longitudinal vibration problem. Therefore, we proceed with more general methods. We first recast our equation of motion into a state space representation of the first-order form by defining a state and excitation vector

$$\mathbf{w}(x,t) = \begin{pmatrix} u_t(x,t) \\ u_x(x,t) \end{pmatrix}, \quad \mathbf{q}(x,t) = \begin{pmatrix} p(x,t) \\ 0 \end{pmatrix} \tag{11}$$

and a matrix differential operator

$$\mathbf{T} = \begin{bmatrix} 0 & c^2 \frac{\partial}{\partial x} \\ \frac{\partial}{\partial x} & 0 \end{bmatrix}. \tag{12}$$

With these definitions our equation of motion attains the first-order form

$$\mathbf{w}_t(x,t) = \mathbf{T}\mathbf{w}(x,t) + \mathbf{q}(x,t) \tag{13}$$

with initial

$$\mathbf{w}(x,0) = \begin{pmatrix} g(x) \\ f'(x) \end{pmatrix} \tag{14}$$

and boundary conditions

$$\begin{bmatrix} -h_1/c & 1 \end{bmatrix} \mathbf{w}(0,t) = 0 \tag{15}$$

$$[h_2/c \quad 1] \, \mathbf{w}(L,t) = 0.$$

Assuming a separable solution

$$\mathbf{w}(x,t) = \mathbf{u}(x)e^{\lambda t} \tag{16}$$

where $\mathbf{u}(x) = \begin{pmatrix} u_1(x) \\ u_2(x) \end{pmatrix}$, Eq. (13) leads to BVP of the form

$$\mathbf{T}\mathbf{u}(x) = \lambda \mathbf{u}(x)$$

$$u_2(0) = \frac{h_1}{c} u_1(0) \tag{17}$$

$$u_2(L) = -\frac{h_2}{c} u_1(L)$$

where $\lambda$ is an eigenvalue while $\mathbf{u}(x)$ is an eigenvector associated with this eigenvalue. Comparing this BVP to the one in Eq. (5), we see that the eigenvalue is not contained in the boundary conditions. It will now be possible to obtain the orthogonality condition which decouples the state space representation Eq. (13).

In view of the above we define a vector space

$$\mathbf{X} \equiv \left\{ \mathbf{u}(x) \; : \; \mathbf{u}(x) = \begin{pmatrix} u_1(x) \\ u_2(x) \end{pmatrix}, \; u_1, u_2 : [0,L] \to \mathbb{C} \; \& \; x \in \mathbb{R} \right\}$$

where $\mathbb{C}$ and $\mathbb{R}$ stand for the set of real and complex numbers respectively. Therefore vector space $\mathbf{X}$ is a 2-tuple of complex valued functions of real variable $x$ defined on the interval $[0, L]$.

We define an inner product $<.,.>: \mathbf{X} \times \mathbf{X} \to F$ as

$$<\mathbf{u}(x), \mathbf{v}(x)> = \int_0^L \left[ \bar{u}_1(x) v_1(x) + \bar{u}_2(x) v_2(x) \right] dx \tag{18}$$

where the bars over $u_1$ and $u_2$ represent the complex conjugation and $F$ stands for the scalar field. It can be easily verified that $<\mathbf{u}(x), \mathbf{v}(x)>$ satisfies positive-definiteness, linearity in the second argument and conjugate-symmetry which are the properties needed for the definition of an inner product. The inner product $<\mathbf{u}(x), \mathbf{v}(x)>$ applied to elements of $\mathbf{X}$ forms an inner product space.

We now proceed to determine the adjoint of operator $\mathbf{T}$ denoted by $\mathbf{T}^*$ where $\mathbf{T}$ and $\mathbf{T}^*$ are defined over some domains $\mathcal{D}(\mathbf{T})$ (or just $\mathcal{D}$) and $\mathcal{D}^*(\mathbf{T}^*)$ (or just $\mathcal{D}^*$) respectively, see §1.6 in [7]. More specifically we look for a vector $\mathbf{g}(x)$ such that

$$<\mathbf{g}(x), \mathbf{u}(x)> \; = \; <\mathbf{v}(x), \mathbf{T}\mathbf{u}(x)> \tag{19}$$

where $\mathbf{u}$ is any element of $\mathcal{D}$, $\mathbf{v}$ is any element of $\mathcal{D}^*$ and $\mathbf{g}(x) = \begin{pmatrix} g_1(x) \\ g_2(x) \end{pmatrix}$. Applying the definition of the inner product to Eq. (19) yields

$$\int_0^L \left[ \bar{g}_1(x)u_1(x) + \bar{g}_2(x)u_2(x) \right] dx = \int_0^L \left[ \bar{v}_1(x) c^2 u_2'(x) + \bar{v}_2(x) u_1'(x) \right] dx \qquad (20)$$

$$= \left[ \bar{v}_1(x) c^2 u_2(x) + \bar{v}_2(x) u_1(x) \right]_0^L - \int_0^L \left[ \bar{v}_1'(x) c^2 u_2(x) + \bar{v}_2'(x) u_1(x) \right] dx.$$

Comparing the integral on the left-hand side of Eq. (20) with the integral on the right-hand side, we see that in order for the equation to hold, it must be that $\bar{g}_1(x) = -\bar{v}_2'(x)$ and $\bar{g}_2(x) = -c^2 \bar{v}_1'(x)$. Furthermore, the first term, the boundary contributions, or so-called "surface term", on the right-hand side of Eq. (20) must be equal to zero. This leads to a definition of $\mathbf{T}^*$ as

$$\mathbf{T}^* = -\begin{bmatrix} 0 & \dfrac{\partial}{\partial x} \\ c^2 \dfrac{\partial}{\partial x} & 0 \end{bmatrix}. \qquad (21)$$

We obtain the boundary conditions for the adjoint problem by setting the surface term to zero

$$\bar{v}_1(L) c^2 u_2(L) + \bar{v}_2(L) u_1(L) - \bar{v}_1(0) c^2 u_2(0) - \bar{v}_2(0) u_1(0) = 0. \qquad (22)$$

Substituting the boundary conditions from Eq. (17) into Eq. (22) and requiring that Eq.(22) is satisfied for all possible $u_1(L)$ and $u_1(0)$, we obtain the boundary conditions for the adjoint eigenvalue problem as $-c h_2 \bar{v}_1(0) + \bar{v}_2(0) = 0$ and $c h_1 \bar{v}_1(0) + \bar{v}_2(0) = 0$. The adjoint eigenvalue problem then becomes

$$\mathbf{T}^* \mathbf{v}(x) = \mu \, \mathbf{v}(x) \qquad (23)$$

$$v_2(0) = -c h_1 v_1(0)$$

$$v_2(L) = c h_2 v_1(L)$$

where $\mu$ is the eigenvalue of operator $\mathbf{T}^*$. We immediately see that $\mathbf{T}^* = -\mathbf{T}^T$ implying $\mathbf{T} \neq \mathbf{T}^*$ and demonstrating that $\mathbf{T}$ is a non-self-adjoint operator. This establishes the non-self-adjointness of BVP, Eq. (5).

We thus obtain a relation (24) which is satisfied when $\mathbf{u}$ is any element of $\mathcal{D}$ and $\mathbf{v}$ is any element of $\mathcal{D}^*$ where $\mathcal{D}$ and $\mathcal{D}^*$ represent the set of eigenvectors satisfying associated boundary conditions for each operator $\mathbf{T}$ and $\mathbf{T}^*$ respectively.

$$< \mathbf{v}(x), \mathbf{T}\mathbf{u}(x) > - < \mathbf{T}^* \mathbf{v}(x), \mathbf{u}(x) > = 0 \qquad (24)$$

We now proceed with obtaining the required orthogonality condition. It can be shown that an eigenvector of $\mathbf{T}$ corresponding to the eigenvalue $\lambda$ is orthogonal to every eigenvector of $\mathbf{T}^*$ corresponding to the eigenvalue $\mu$ for which $\mu \neq \bar{\lambda}$ and that the set of $\lambda$'s of operator $\mathbf{T}$ is the same as the set of $\mu$'s of operator $\mathbf{T}^*$ (see §2.4 in [7]).

To determine eigenvectors of $\mathbf{T}$ we revert to BVP specified in Eq. (5). The eigenvalues of that problem are the same as eigenvalues of the problem specified in Eq. (17) excluding the eigenvalue of zero. This is because $\lambda = 0$ for BVP in Eq. (17) has only a trivial solution which cannot be an eigenvector. Eigenvectors for both operators $\mathbf{T}$ and $\mathbf{T}^*$ can now readily be obtained by using Eq. (7) in view of Eq. (17) as

$$\mathbf{u}_n(x) = \begin{pmatrix} \cosh\left(\frac{\lambda_n}{c}x\right) + h_1 \sinh\left(\frac{\lambda_n}{c}x\right) \\ \frac{1}{c}\left[\sinh\left(\frac{\lambda_n}{c}x\right) + h_1 \cosh\left(\frac{\lambda_n}{c}x\right)\right] \end{pmatrix} = \begin{pmatrix} u_{1,n}(x) \\ u_{2,n}(x) \end{pmatrix}, \tag{25}$$

and in view of Eqs. (23) and (25) as

$$\mathbf{v}_n(x) = \begin{pmatrix} \cosh\left(\frac{\overline{\lambda}_n}{c}x\right) + h_1 \sinh\left(\frac{\overline{\lambda}_n}{c}x\right) \\ -c\left[\sinh\left(\frac{\overline{\lambda}_n}{c}x\right) + h_1 \cosh\left(\frac{\overline{\lambda}_n}{c}x\right)\right] \end{pmatrix} = \begin{pmatrix} \overline{u}_{1,n}(x) \\ -c^2 \overline{u}_{2,n}(x) \end{pmatrix}. \tag{26}$$

Vectors $\mathbf{u}_n$ and $\mathbf{v}_n$ represent infinite orthogonal bases in $\mathbf{X}$ which can be used to expand any two-component vector $\mathbf{F}$ into a series

$$\mathbf{F}(x) = \sum_{n=-\infty}^{\infty} \alpha_n \mathbf{u}_n(x) \tag{27}$$

where

$$\alpha_n = \frac{<\mathbf{v}_n(x), \mathbf{F}(x)>}{L(1-h_1^2)}. \tag{28}$$

The denominator in Eq. (28) is the result of evaluating $<\mathbf{v}_n(x), \mathbf{u}_n(x)>$.

Eq. (27) represents a generalized Fourier series for which it can be shown that the coefficients $\alpha_n$ are the best possible choices in the mean square sense. Furthermore, they are components of the projections onto the basis $\mathbf{u}_n$ and since the projections are never greater than the original, no piece of the Fourier series can have the square norm greater than the norm of $\mathbf{F}$. In the limit as $n \to \pm\infty$ the series on the left-hand side of Eq. (27) reconstructs $\mathbf{F}$ demonstrating the completeness of the basis. The proof of this can be found in [8] (Theorem 2.1, condition 2.11) as well as in [5].

## 5. Modal analysis

We proceed by using the results of the previous section and decoupling the system of equations given in Eq. (13). Assuming a separable solution

$$\mathbf{w}(x,t) = \sum_{r=-\infty}^{\infty} \eta_r(t) \mathbf{u}_r(x) \tag{29}$$

and substituting it into the system (13) yields

$$\sum_{r=-\infty}^{\infty} \dot{\eta}_r(t) \mathbf{u}_r(x) = \mathbf{T}\left(\sum_{r=-\infty}^{\infty} \eta_r(t) \mathbf{u}_r(x)\right) + \mathbf{q}(x,t) \qquad (30)$$

$$= \sum_{r=-\infty}^{\infty} \eta_r(t) \lambda_r \mathbf{u}_r(x) + \mathbf{q}(x,t).$$

We now take the inner product of the left-hand and right-hand side of the equation with $<\mathbf{v}_q(x), \cdot>$ to obtain a decoupled set of equation. Reverting back to index $r$ from $q$, we write a decoupled set of equations as

$$\dot{\eta}_r(t) = \eta_r(t) \lambda_r + \frac{1}{(1-h_1^2)L} \int_0^L p(\xi,t) u_{1,r}(\xi) d\xi, \qquad \forall r \qquad (31)$$

where $u_{1,r}(x)$ is the first component of $\mathbf{u}_r(x)$ vector. Applying the Laplace transform to Eq. (31) we obtain

$$\tilde{\eta}_r(s) = \frac{\eta_r(0)}{s-\lambda_r} + \frac{1}{(1-h_1^2)L\,(s-\lambda_r)} \int_0^L \tilde{p}(\xi,s) u_{1,r}(\xi) d\xi, \qquad (32)$$

where $\tilde{\eta}_r$ and $\tilde{p}$ represent Laplace transforms of time-dependent variables with $s$ being a parameter.

On the other hand, we know that the Laplace transform assumes expansion

$$\mathcal{L}\{\mathbf{w}(x,t)\} = \mathcal{L}\left\{\begin{pmatrix} u_t(x,t) \\ u_x(x,t) \end{pmatrix}\right\} = \sum_{r=-\infty}^{\infty} \mathcal{L}\{\eta_r(t)\} \mathbf{u}(x). \qquad (33)$$

We obtain

$$\begin{pmatrix} s\tilde{u}(x,s) - u(x,0) \\ \tilde{u}_x(x,s) \end{pmatrix} = \sum_{r=-\infty}^{\infty} \tilde{\eta}_r(s) \begin{pmatrix} u_{1,r}(x) \\ u_{2,r}(x) \end{pmatrix}. \qquad (34)$$

The first component of the vector on both sides of Eq. (34) provides the Laplace transform of the response of the bar as

$$\tilde{u}(x,s) = \frac{u(x,0)}{s} + \sum_{r=-\infty}^{\infty} \frac{\tilde{\eta}_r(s)}{s} u_{1,r}(x). \qquad (35)$$

In view of Eq. (32) we obtain

$$\tilde{u}(x,s) = \frac{u(x,0)}{s} + \sum_{r=-\infty}^{\infty} \frac{\eta_r(0)}{s(s-\lambda_r)} u_{1,r}(x) + \sum_{r=-\infty}^{\infty} \left(\frac{1}{(1-h_1^2)Ls(s-\lambda_r)} \int_0^L \tilde{p}(\xi,s) u_{1,r}(\xi) d\xi\right) u_{1,r}(x). \qquad (36)$$

The response in the time domain is the inverse Laplace transform of Eq. (36)

$$u(x,t) = u(x,0) + \sum_{r=-\infty}^{\infty} \eta_r(0) \left[ -\frac{1}{\lambda_r}(1 - e^{\lambda_r t}) \right] u_{1,r}(x) \qquad (37)$$

$$+ \frac{1}{(1-h_1^2)L} \sum_{r=-\infty}^{\infty} \left( \int_0^t \left[ -\frac{1}{\lambda_r}(1 - e^{\lambda_r(t-\tau)}) \right] \int_0^L p(\xi,\tau) u_{1,r}(\xi) d\xi \, d\tau \right) u_{1,r}(x)$$

where $\eta_r(0)$ can be calculated from the initial conditions

$$\mathbf{w}(x,0) = \begin{pmatrix} u_t(x,0) \\ u_x(x,0) \end{pmatrix} = \begin{pmatrix} g(x) \\ f'(x) \end{pmatrix} = \sum_{r=-\infty}^{\infty} \eta_r(0) \begin{pmatrix} u_{1,r}(x) \\ u_{2,r}(x) \end{pmatrix}. \qquad (38)$$

Applying the inner product $<\mathbf{v}_q(x), .>$ to Eq. (38) and reverting back to index $r$ from $q$, we obtain decoupled coefficients of expansion as

$$\eta_r(0) = \frac{1}{(1-h_1^2)L} \int_0^L \left[ g(\xi) u_{1,r}(\xi) - c^2 f'(\xi) u_{2,r}(\xi) \right] d\xi, \quad \forall r. \qquad (39)$$

Simplifying the integral by partial integration and using boundary conditions from Eq. (17), we obtain $\eta_r(0)$ written in terms of $u_{1,r}$ as

$$\eta_r(0) = \frac{1}{(1-h_1^2)L} \left\{ \int_0^L [\lambda_r f(\xi) + g(\xi)] u_{1,r}(\xi) d\xi + c f(L) h_2 u_{1,r}(L) + c f(0) h_1 u_{1,r}(0) \right\}, \quad \forall r. \qquad (40)$$

We now substitute Eq. (40) into Eq. (37) and, in view of Appendix and continuity of $u_{1,r}(\xi)$ and $p(\xi,\tau)$ on $[0,L]$ together with interchanging the summation and integral signs, we calculate individual terms as follows:

$$f(x) + \frac{1}{L(1-h_1^2)} \sum_{r=-\infty}^{\infty} \int_0^L \lambda_r f(\xi) u_{1,r}(\xi) d\xi \left( -\frac{1}{\lambda_r} \right) u_{1,r}(x) = 0, \qquad (41)$$

$$\frac{1}{(1-h_1^2)L} \sum_{r=-\infty}^{\infty} \left( \left( -\frac{1}{\lambda_r} \right) \int_0^L g(\xi) u_{1,r}(\xi) d\xi \right) u_{1,r}(x) = \qquad (42)$$

$$\int_0^L g(\xi) \sum_{r=-\infty}^{\infty} \left( -\frac{u_{1,r}(\xi) u_{1,r}(x)}{\lambda_r (1-h_1^2) L} \right) d\xi = \frac{1}{c(h_1+h_2)} \int_0^L g(\xi) d\xi,$$

$$\frac{1}{(1-h_1^2)L} \sum_{r=-\infty}^{\infty} \left( \left( -\frac{1}{\lambda_r} \right) c f(L) h_2 u_{1,r}(L) \right) u_{1,r}(x) = c f(L) h_2 \sum_{r=-\infty}^{\infty} \left( -\frac{u_{1,r}(L) u_{1,r}(x)}{\lambda_r (1-h_1^2) L} \right) = f(L) \frac{h_2}{h_1+h_2}, \qquad (43)$$

$$\frac{1}{(1-h_1^2)L}\sum_{r=-\infty}^{\infty}\left(\int_0^t\left(-\frac{1}{\lambda_r}\right)\int_0^L p(\xi,\tau)u_{1,r}(\xi)d\xi\,d\tau\right)u_{1,r}(x)= \tag{44}$$

$$\int_0^t\int_0^L p(\xi,\tau)\sum_{r=-\infty}^{\infty}\left(-\frac{u_{1,r}(\xi)u_{1,r}(x)}{\lambda_r(1-h_1^2)L}\right)d\xi\,d\tau = \frac{1}{c(h_1+h_2)}\int_0^t\int_0^L p(\xi,\tau)d\xi\,d\tau.$$

Finally, in view of Eqs. (41) – (44) we can write the vibratory response as

$$u(x,t)=\frac{1}{c(h_1+h_2)}\left[\int_0^L g(\xi)d\xi+c\,f(L)h_2\right]+\frac{c\,f(L)h_2}{(1-h_1^2)L}\sum_{r=-\infty}^{\infty}\frac{u_{1,r}(L)u_{1,r}(x)}{\lambda_r}e^{\lambda_r t} \tag{45}$$

$$+\frac{1}{(1-h_1^2)L}\sum_{r=-\infty}^{\infty}\frac{u_{1,r}(x)}{\lambda_r}e^{\lambda_r t}\int_0^L u_{1,r}(\xi)[\lambda_r f(\xi)+g(\xi)]d\xi$$

$$+\frac{1}{c(h_1+h_2)}\int_0^t\int_0^L p(\xi,\tau)d\xi\,d\tau+\frac{1}{(1-h_1^2)L}\sum_{r=-\infty}^{\infty}\frac{u_{1,r}(x)}{\lambda_r}\int_0^t e^{\lambda_r(t-\tau)}\int_0^L p(\xi,\tau)u_{1,r}(\xi)d\xi\,d\tau$$

$$h_1+h_2\neq 0,\ h_1\neq\pm 1.$$

It should be pointed out that we are, of course, only interested in the real part of Eq. (45). However, there is no particular need to extract only the real part from the response since the imaginary part must necessarily be zero due to all of the initial conditions being real as well as the force applied.

## 6. Verifications

In this section we investigate the general response further and demonstrate that it contains the solutions to the other three cases: free-free, fixed-damper, fixed-fixed and fixed-free bar. Free-free, fixed-fixed, fixed-free bar are textbook cases, while fixed-damper bar case is a case with non-self-adjoint boundary conditions. However, we first verify that the obtained general response is self-consistent. In other words, we expect that providing an initial force impulse along the length of the bar at rest is equivalent to imposing initial velocity of the bar without any force. The vibratory response of the bar in both cases should be the same provided that $f(x)$ is zero. To this end, we first determine what the initial velocity of the bar is to an impulsive force $p(x,t)=p_0\delta(t)$ where $p_0$ is a constant and $\delta(t)$ is Dirac's delta function. Using the impulse-momentum equation yields

$$\int_0^{0^+}\rho A_0 u_{tt}dt=\int_0^{0^+}\rho A_0 p_0\delta(t)dt\ \Rightarrow\ u_t(x,0^+)=p_0. \tag{46}$$

Therefore, imposing impulsive force $p(x,t)=p_0\delta(t)$ must be equivalent to imposing initial velocity $u_t(x,0^+)=p_0$. To confirm this we substitute $g(x)=p_0$, $p(x,t)=0$ and $f(x)=0$ into Eq. (45), and calculate the response of the bar

$$u(x,t)=\frac{p_0 L}{c(h_1+h_2)}+\frac{1}{(1-h_1^2)L}\sum_{r=-\infty}^{\infty}\frac{p_0 u_{1,r}(x)}{\lambda_r}e^{\lambda_r t}\int_0^L u_{1,r}(x)d\xi. \tag{47}$$

On the other hand, substituting $g(x)=0$, $p(x,t)=p_0\delta(t)$ and $f(x)=0$ yields the response

$$u(x,t) = \frac{1}{c(h_1+h_2)} \int_0^t \int_0^L p_0 \delta(\tau) d\xi d\tau + \frac{1}{(1-h_1^2)L} \sum_{r=-\infty}^{\infty} \frac{u_{1,r}(x)}{\lambda_r} \int_0^t e^{\lambda_r(t-\tau)} \int_0^L p_0 \delta(\tau) u_{1,r}(\xi) d\xi d\tau \quad (48)$$

$$= \frac{p_0}{c(h_1+h_2)} \int_0^L d\xi \int_0^t \delta(\tau) d\tau + \frac{1}{(1-h_1^2)L} \sum_{r=-\infty}^{\infty} \frac{u_{1,r}(x)}{\lambda_r} e^{\lambda_r t} \int_0^t e^{-\lambda_r \tau} \delta(\tau) d\tau \int_0^L p_0 u_{1,r}(\xi) d\xi$$

As expected, these two responses are the same since the integrals of $\delta(\tau)$ evaluate to unity. We further verify that Eq. (45) is correctly predicting the motion of a rigid bar. As the stiffness of the continuous bar increases, the vibrations should reduce and the solution should converge to that of a dashpot-mass system. We first note that since $h_i = c_i \sqrt{1/E\rho}/A_0$ we have $h_i < 1$ when $E > c_i^2/(\rho A_0^2)$ where $i = 1, 2$. Therefore, as $E \to +\infty$ for $E > c_i^2/(\rho A_0^2)$ we have $Arg\left((1-h_1)(1-h_2)/((1+h_1)(1+h_2))\right) = 0$ and consequently, the imaginary part of eigenvalues in Eq. (10) becomes $n\pi c/L$. Next, we observe that by Eqs. (10) and (7) for the real part of the eigenvalues we have $\lim_{E \to +\infty} \text{Re}(\lambda_r) = -(c_1 + c_2)/(\rho AL)$ as well as $\lim_{E \to +\infty} u_{1,0}(x) = 1$ and $\lim_{E \to +\infty} u_{1,r}(x) = \cos(r\pi x/L)$. We are now ready to evaluate the limit of the right hand side of Eq. (45) as $E \to +\infty$ given no initial displacement, initial velocity $g(x) = \dot{x}_0 = const$ and no external force on the continuous bar. This results in

$$\lim_{E \to +\infty} \left\{ \frac{1}{c(h_1+h_2)} \left[ \int_0^L g(\xi) d\xi \right] + \frac{u_{1,0}(x) e^{\lambda_0 t}}{(1-h_1^2) L \lambda_0} \int_0^L u_{1,0}(\xi) g(\xi) d\xi \right. \quad (49)$$

$$\left. + \frac{1}{(1-h_1^2) L} \sum_{r=1}^{\infty} \left[ \frac{u_{1,r}(x)}{\lambda_r} e^{\lambda_r t} \int_0^L u_{1,r}(\xi) g(\xi) d\xi + \frac{u_{1,-r}(x)}{\lambda_{-r}} e^{\lambda_{-r} t} \int_0^L u_{1,-r}(\xi) g(\xi) d\xi \right] \right\}$$

$$= \frac{\rho AL \dot{x}_0}{c_1 + c_2} \left( 1 - e^{-\frac{c_1+c_2}{\rho AL} t} \right) + \frac{\dot{x}_0}{\pi} \sum_{r=1}^{\infty} \frac{1}{r} \sin(r\pi) \cos\left(\frac{r\pi x}{L}\right) \lim_{E \to +\infty} \left[ \frac{e^{\lambda_r t}}{\lambda_r} + \frac{e^{\lambda_{-r} t}}{\lambda_{-r}} \right].$$

The limit on the right hand side of Eq. (49) becomes

$$\lim_{E \to +\infty} \left[ \frac{e^{\lambda_r t}}{\lambda_r} + \frac{e^{\lambda_{-r} t}}{\lambda_{-r}} \right] \leq 2 e^{-\frac{c_1+c_2}{\rho AL} t} \lim_{E \to +\infty} \frac{1}{|\lambda_r|} \left[ \cos(\text{Im}(\lambda_r)t) + \sin(\text{Im}(\lambda_r)t) \right] = 0 \quad (50)$$

where $\text{Im}(\lambda_r)$ represents the imaginary part of eigenvalues. Therefore, the second term on the right hand side of Eq. (49) vanishes while the first term is exactly the solution to the motion of a rigid uniform bar with mass $\rho AL$ and two dampers attached. We further observe that as $E \to +\infty$ Eq. (49) indicates that all the eigenmodes vanish except the rigid body mode and consequently the rigid body motion of the bar is confirmed.

We now demonstrate that the eigenfunctions for free-free, fixed-damper, fixed-fixed and fixed-free bar can be obtained by considering the vibratory response of damper-damper bar.

To obtain the eigenfunction for free-free bar we simply let $h_1 \to 0$ and $h_2 \to 0$ in Eq. (10) and we choose the appropriate imaginary part. This provides us with the limit for $\lambda_n$ as $\lim\limits_{\substack{h_1 \to 0 \\ h_2 \to 0}} \lambda_n = \frac{cn\pi}{L} i$.

Substitution this into Eq. (7) and taking the limit we obtain $\lim\limits_{\substack{h_1 \to 0 \\ h_2 \to 0}} \phi_n = B\cosh\left(\frac{1}{c}\frac{cn\pi x}{L}i\right) = B\cos\left(\frac{n\pi x}{L}\right)$.

Letting $h_1 \to 0$ only yields the case discussed in [5].

We see immediately that we obtain a well known eigenfunction for free-free bar which is real since free-free bar has the self-adjoint boundary conditions and therefore $n = 1,2,3,...$ Since this case has non-restraining boundaries as well, there must be a rigid body mode which corresponds to $\phi_r = A_r = const$.

The eigenfunction for the fixed-damper bar case is obtainable by renormalizing Eq. (7) by $h_1$ and letting $h_1 \to +\infty$. As another test, it can also be inferred by examining any summation term in the vibratory response. We first determine the eigenvalues for this case by letting $h_1 \to +\infty$ in Eq. (10). This yields

$$\lambda_r = \frac{c}{2L}\left[\ln\left|\frac{(1-h_2)}{(1+h_2)}\right| + i\left(Arg\left(-\frac{(1-h_2)}{(1+h_2)}\right) + 2r\pi\right)\right], \quad r \in \mathbb{Z}. \tag{51}$$

Now we take one of the products of eigenfunction under the summation and integral signs from the vibratory response and let $h_1 \to +\infty$ to obtain

$$\lim_{h_1 \to +\infty} \frac{u_{1,r}(x)u_{1,r}(\xi)}{(1-h_1^2)L} = \lim_{h_1 \to +\infty} \frac{\left[\cosh\left(\frac{\lambda_r}{c}x\right) + h_1\sinh\left(\frac{\lambda_r}{c}x\right)\right]\left[\cosh\left(\frac{\lambda_r}{c}\xi\right) + h_1\sinh\left(\frac{\lambda_r}{c}\xi\right)\right]}{(1-h_1^2)L} \tag{52}$$

$$= -\frac{\sinh\left(\frac{\lambda_r}{c}x\right)\sinh\left(\frac{\lambda_r}{c}\xi\right)}{L}.$$

We immediately see that the eigenfunction for this case must be $\sinh\left(\frac{\lambda_r}{c}x\right)$ which is complex valued and since this case represents a non-self-adjoint problem $r \in \mathbb{Z}$. However, because one boundary is restrained, there is no rigid body mode present and therefore there is no need for a constant eigenfunction. From this case it is trivial to arrive to the case of fixed-fixed bar by taking a limit $h_2 \to \infty$ in Eqs. (51) and then $\lim\limits_{h_2 \to +\infty} \sinh\left(\frac{\lambda_r}{c}x\right)$.

Finally, we determine the eigenfunction for fixed-free bar from the previous case by letting $h_2 \to 0$ in Eq. (51) and substituting the resulting eigenvalue into $\sinh\left(\frac{\lambda_r}{c}x\right)$. Therefore we obtain the eigenvalue $\lim_{h_2 \to 0} \lambda_r = \frac{c}{L}\frac{(2r+1)\pi}{2}i$ and the eigenfunction

$$\sinh\left(\frac{1}{c}\frac{c}{L}\frac{(2r+1)\pi}{2}ix\right) = i\sin\left(\frac{(2r+1)\pi}{2L}x\right). \tag{53}$$

From Eq. (53) we see immediately that, since $i$ is just a multiplicative constant, the eigenfunction is real imposing $r = 0,1,2,3,...$ Thus, we obtained the well known eigenfunction for the fixed-free bar as expected.

Any of the above eigenfunctions can be used in the general response to determine the response for the desired case provided that proper limits and boundary conditions are used in Eq. (45). To demonstrate that, and as a final verification for this section, we determine the frequency response for fixed-damper bar by applying a harmonic force $\frac{F_0}{\rho A_0}e^{i\omega t}\delta(x - x_f)$ at some point $x_f$ where $F_0$ is a constant force per unit length as given in [4]. The vibratory response for this case is

$$u(x,t) = \frac{1}{c(h_1+h_2)}\int_0^t\int_0^L \frac{F_0}{\rho A_0}e^{i\omega\tau}\delta(\xi - x_f)d\xi\, d\tau \tag{54}$$

$$+ \frac{1}{(1-h_1^2)L}\sum_{r=-\infty}^{\infty}\frac{u_{1,r}(x)}{\lambda_r}\int_0^t e^{\lambda_r(t-\tau)}\int_0^L \frac{F_0}{\rho A_0}e^{i\omega\tau}\delta(\xi - x_f)u_{1,r}(\xi)d\xi\, d\tau$$

$$= \frac{1}{c(h_1+h_2)}\int_0^t\int_0^L \frac{F_0}{\rho A_0}e^{i\omega\tau}\delta(\xi - x_f)d\xi\, d\tau$$

$$+ \frac{1}{(1-h_1^2)L}\sum_{r=-\infty}^{\infty}\frac{u_{1,r}(x)}{\lambda_r}e^{\lambda_r t}\frac{F_0}{\rho A_0}u_{1,r}(x_f)\left[\frac{e^{-\tau(\lambda_r - i\omega)}}{i\omega - \lambda_r}\right]_0^t.$$

We already determined the product of eigenfunctions from Eq. (52) when $h_1 \to +\infty$ and therefore we can immediately write the limit of Eq. (54) in view of eigenvalues given in Eq. (51). This yields

$$u(x,t) = -\frac{F_0}{\rho A_0 L}\sum_{r=-\infty}^{\infty}\frac{\sinh\left(\frac{\lambda_r x}{c}\right)\sinh\left(\frac{\lambda_r x_f}{c}\right)e^{i\omega t}}{\lambda_r(i\omega - \lambda_r)}. \tag{55}$$

Here, only the steady state exponential terms are retained. This result is the same as reported in [4].

## 7. Numerical results

In this section we investigate numerically the general response provided in Eq. (45). First, we are interested to know how well the eigenfunctions approximate arbitrary functions on [0, $L$]. To determine this we assume a closed-form expression

$$u_e(x,t) = x^2\left[1 - \frac{(h_2 L - 2c)x}{L(h_2 L - 3c)}\right]e^{-t} \tag{56}$$

which satisfies the boundary conditions and produces the required forcing function

$$p_e(x,t) = \frac{2e^{-t}}{L(3c - h_2 L)}\left[\left(c^2 - \frac{1}{2}x^2\right)h_2 L^2 + \left(\frac{3}{2}x^2 c + \frac{1}{2}x^3 h_2 - 3c^3 - 3c^2 x h_2\right)L - x^3 c + 6c^3 x\right]. \tag{57}$$

Using the response Eq. (45), the displacement error is given with $\varepsilon(x,t) = u_e(x,t) - u(x,t)$ Figure 2 depicts the error $\varepsilon(x, 0.3)$ for which the sums in Eq. (45) are calculated as $\sum_{r=-k}^{k}$. The error decreases rapidly as $k$ increases and is within 0.0005 for $k = 15$.

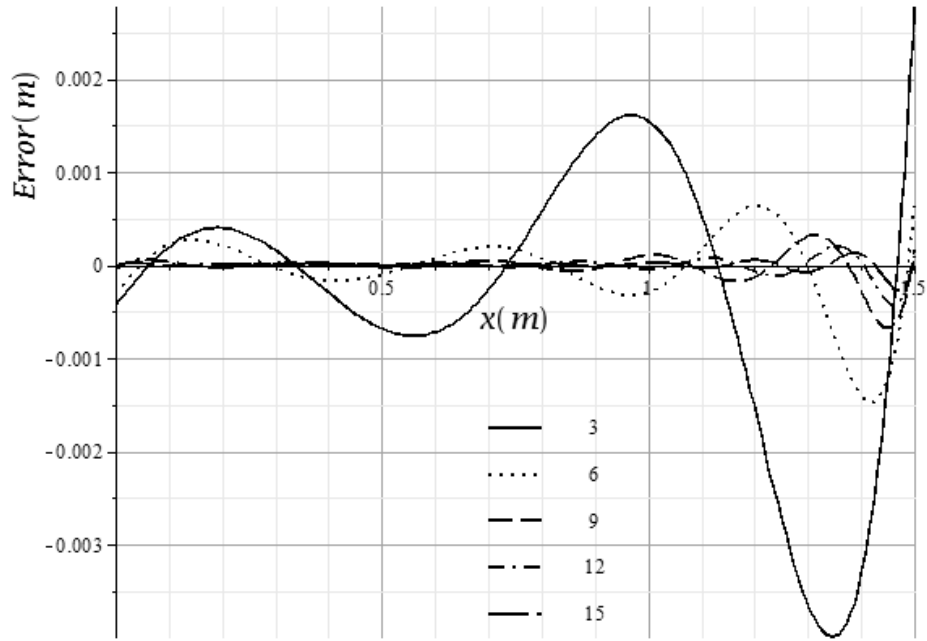

**Figure 2. Error at $t = 0.3$ for $h_1 = 0.3$, $h_2 = 0.7$, $c = 1.8$, $L = 1.5$ and $k = 3, 6, 9, 12, 15$**

In Figure **3** a response for the same parameters as in Figure 2 until time $t = 2s$ is depicted. It can be observed that the constrains acting on the bar are such that its left end does not move for all time indicating that the damper attached at $x = 0$ is not seeing any force. Therefore, it can be concluded that the series solution Eq. (45), must yield the same displacement of the bar at any particular instant of time for any value of $h_1$ that is not special (section 3).

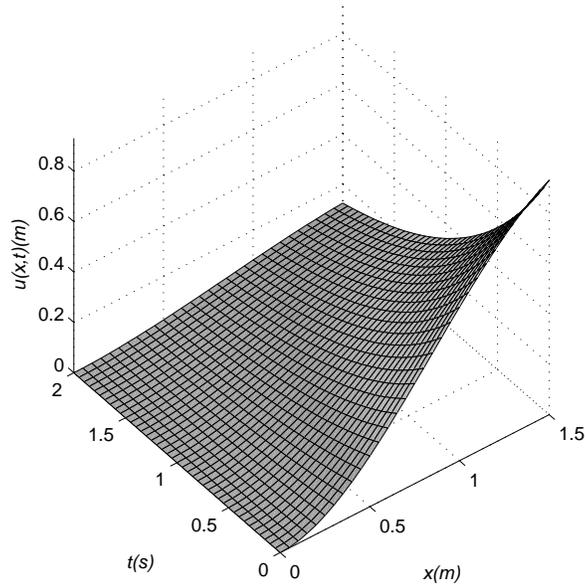

**Figure** 3. **The response of the system corresponding to imposed** $p_e(x,t)$ **and initial conditions obtained from** $u_e(x,t)$ **up to time** *t* = **2s.**

We now wish to see how the eigenfunction expansion fares when compared with FEM numerical simulations. For FEM we use linear shape functions $[1-x/l \quad x/l]$ where $l$ is the length of the simple linear spring finite element. We determine energy equivalent nodal forces from $\rho A_0 \int_0^l \begin{bmatrix} 1-x/l \\ x/l \end{bmatrix} p(x+x_i,t)dx$ where $x_i$ is the nodal distance from the left end of the bar. The assembled system of equations, after simplifications, can be written as

$$\frac{L}{6nc}\begin{bmatrix} 2 & 1 & 0 & . & . & 0 \\ 1 & 4 & 1 & 0 & . & 0 \\ . & . & . & . & . & . \\ . & . & . & . & . & . \\ 0 & . & 0 & 1 & 4 & 1 \\ 0 & . & . & 0 & 1 & 2 \end{bmatrix}\ddot{U}(t) + \begin{bmatrix} h_1 & 0 & . & . & . & 0 \\ 0 & 0 & . & . & . & 0 \\ . & . & . & . & . & . \\ . & . & . & . & . & . \\ 0 & . & . & . & 0 & 0 \\ 0 & . & . & . & 0 & h_2 \end{bmatrix}\dot{U}(t)$$

$$+\frac{nc}{L}\begin{bmatrix} 1 & -1 & . & . & . & 0 \\ -1 & 2 & -1 & . & . & 0 \\ . & . & . & . & . & . \\ . & . & . & . & . & . \\ 0 & . & . & -1 & 2 & -1 \\ 0 & . & . & 0 & -1 & 1 \end{bmatrix}U(t) = \frac{1}{c}\begin{bmatrix} \int_0^l \left(1-\frac{x}{l}\right)p(x,t)dx \\ \int_0^l \left[\left(1-\frac{x}{l}\right)p(x+x_2,t)+\frac{x}{l}p(x,t)\right]dx \\ . \\ \int_0^l \left[\left(1-\frac{x}{l}\right)p(x+x_i,t)+\frac{x}{l}p(x+x_{i-1},t)\right]dx \\ . \\ \int_0^l \frac{x}{l}p(x+x_n,t)dx \end{bmatrix}$$ (58)

where $n$ is the number of finite elements, $l = L/n$, $i = 1, 2,\ldots, n+1$ and $U(t) = [U_1(t), U_2(t),\ldots, U_{n+1}(t)]^T$ is the displacement vector with displacement components at distances $x_1 = 0, x_2,\ldots, x_{n+1}$ from the left end of the bar. Therefore, the mesh for FEM is simply a division of the full length $L$ of the bar into $n$ linear spring elements of the length $l$.

Figure 4 depicts the response of the system to force $p(x,t) = \sin(6\pi x / L)\sin(\pi t / L)$ up to time 2 s with initial conditions $f(x) = 0.1x(L-x/2)$ and $g(x) = 0$ using the analytical expression and the same set of parameters as in Figure 2.

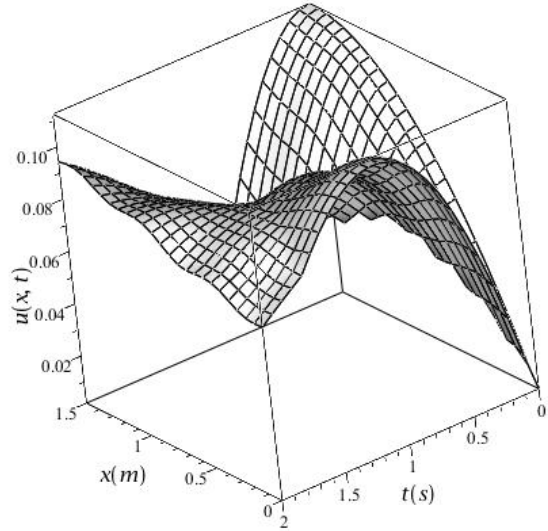

**Figure 4. Response of the system to force** $p(x,t) = \sin(6\pi x/L)\sin(\pi t/L)$ **up to time** *t*=2 s.

Figure 5 shows maximum absolute displacement difference between the analytical method and FEM for the same parameters as in Figure 2 at $t = 0.5$s. FEM was calculated with Eq. (58) with 40 linear spring elements while the analytical response was calculated with Eq. (45) with increasing number of eigenfunctions. Both calculations were performed with the inputs that produced Figure 4. As the number of eigenfunctions increases, the analytical method rapidly converges to the displacement produced by FEM. We find that with just three eigenfunctions the displacement calculated by analytical method is already within 0.011 relative to the displacement of FEM.

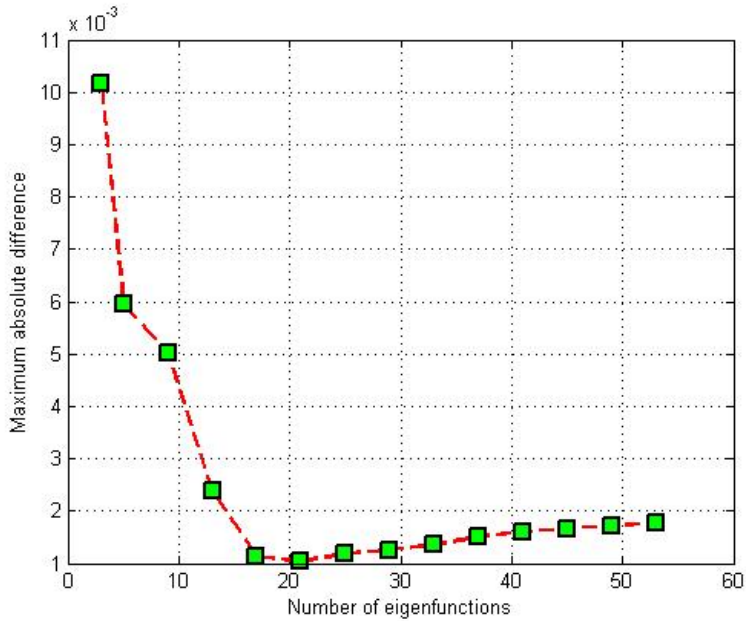

**Figure 5. Maximum absolute difference in displacement between FEM with 40 finite elements and the analytical expression.**

With the further increase of the number of eigenfunctions the displacement difference between the two methods reduces until about 22 eigenfunctions when it starts to increase. The reason for the increase is that with 22 eigenfunctions the analytical method achieves better accuracy than FEM with 40 finite elements; consequently, the displacement difference between the two methods increases. This is further confirmed in Figure 6 where we see that for a more accurate FEM the dip in the absolute displacement difference moved to the right and the difference begins to level off with only a slight increase after about 32 eigenfunctions.

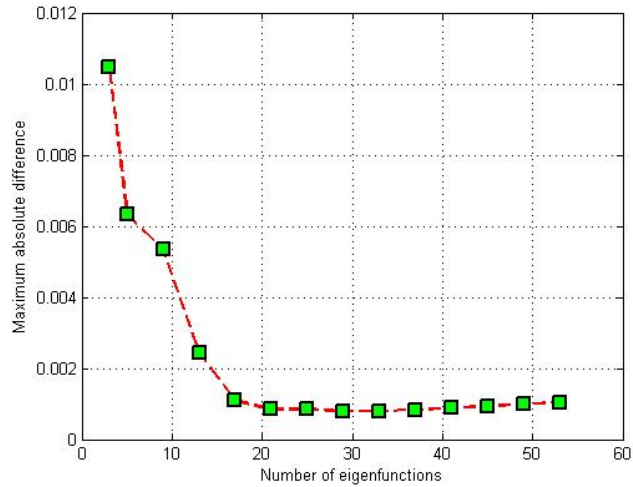

**Figure 6. Maximum absolute difference in displacement between FEM with 60 finite elements and the analytical expression.**

Figure 5 and Figure 6 demonstrate the efficiency of the eigenfunction expansion with using just a few terms in the Fourier expansion. Furthermore, the finite element method performs well, but its disadvantage is that the system has to be simulated from initial to final time unlike the analytical solution which can be computed for any given time. However, to really understand the difference between these two procedures one has to look at the root locus of eigenvalues of FEM in the complex plane as the number of finite elements increases.

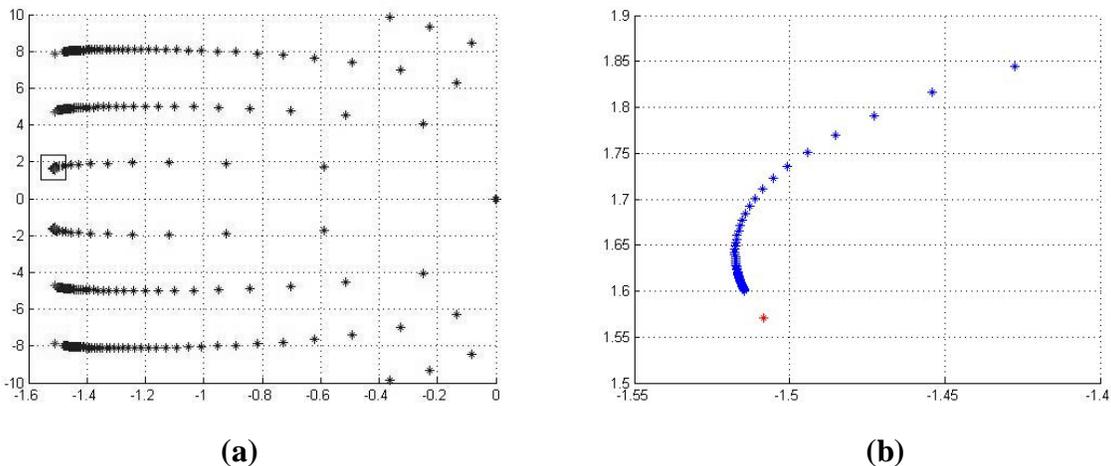

(a)         (b)

**Figure 7. (a) Distribution of eigenvalues in complex plane for FEM calculated with Matlab's eig() function for 1 to 60 finite elements, (b) Magnification around one of the exact eigenvalues.**

Figure 7 shows the part of a root locus in the complex plane for FEM and Eq. (10) for parameters $h_1=.3$, $h_2=1.2$, $c=1$ and $L=1$ for the number of finite elements ranging from 1 to 60. The eigenvalues were calculated with Matlab's eig() function after the system of equations (58) was put into the first order form and the state matrix was obtained from it. The stars situated along the vertical line located on the real axis at –1.508467.on the left side of (a) are the exact eigenvalues For a very small number of elements of FEM, represented with stars on the right side of (a), the root locus of eigenvalues is far from the theoretical ones. As the number of elements increases, the locus of stars is clearly moving towards the theoretical eigenvalues, the left side of (a). The magnification of the root locus in the vicinity of an exact eigenvalue is depicted in (b). Thus, it can be seen that the real parts of eigenvalues of FEM are highly dependent on the number of finite elements and consequently are not constant. This may look unintuitive since our theoretical result states that the real parts of eigenvalues are constant. For a practicing engineer this would have been hard to notice without his or her being aware of the theoretical result in advance. Thus, we demonstrate a need for qualitative understanding of the problem before attempting FEM.

Figure 7 demonstrates how the two methods differ from each other. The Fourier series solution is superior since calculating additional terms in the series amounts to only adding eigenvalues in the complex plane along the vertical line where each eigenvalue has the constant rate of decay. Due to the orthogonality of the basis, the error of approximation reduces in the mean square sense as the new terms are added. This has no damaging effect on the previous calculation. Namely, we do not need to change previously calculated coefficients, but merely to determine additional ones. On the other hand FEM requires a complete recalculation of the eigenvalues and the reduction of error happens at slower rate because the real parts of the eigenvalues of FEM are greater as compared to the real parts of eigenvalues of the analytical response.

As a side observation, in practice, most engineers assume that calculated eigenvalues in FEM are close to the theoretical ones. Although each individual eigenvalue in FEM approaches its theoretical value when the number of elements increases, for any fixed number of elements, no matter how large, the approximate eigenvalues do not appear situated along the vertical line, but rather along what appears to be a parabolic curve. This is because the rate of convergence becomes slower as the imaginary part of eigenvalues increases. Thus, a priori knowledge of the situations involving viscous boundaries is essential for an engineer to properly interpret FEM results.

## 8. Spurious eigenvalues in FEM

We have encountered a phenomenon for which we do not have an explanation at this time. We have observed that for some values of parameters $h_1$ and $h_2$ the stable continuous system becomes unstable when discretized by the FEM method. One set of parameters that produces such a behavior is $h_1 = 0.7$, $h_2 = -1.5$, c =1.5 and $L=1.8$ for which the distribution of eigenvalues is along the line positioned on the real axis at -0.0521513 in addition to one eigenvalue at the origin. It is clear that this continuous system is stable since there are no eigenvalues with positive real parts. If, however, one discretizes the system, an eigenvalue with a positive real part will arise.

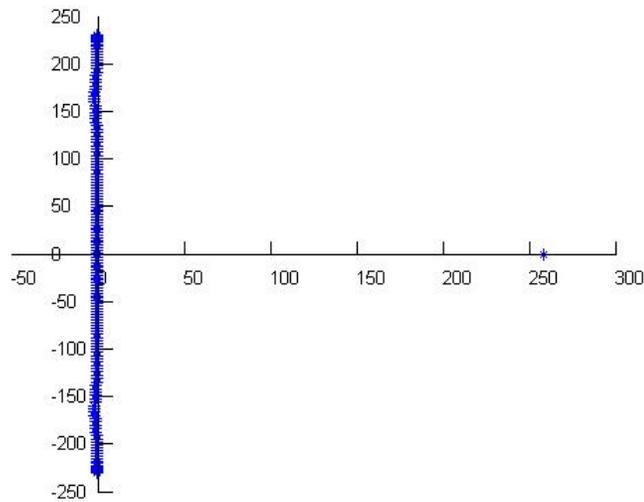

**Figure 8. Eigenvalues from FEM in the complex plane for $h_1 = 0.7$ and $h_2 = -1.5$.**

Figure 8 depicts eigenvalues of the discrete system obtained by FEM. It can be observed that a spurious eigenvalue lies on the positive part of the real axis in the complex plane which forces the discrete system to become unstable. The situation is not improved by increasing the number of finite elements. On the contrary, the real part of the spurious eigenvalue will increase to infinity making the system even more unstable. This is counter-intuitive since we know that the continuous system is equivalent to the discrete one as the number of elements tends to infinity. Therefore, it can be concluded that at least some stability regions for some parameters of the continuous system become so distorted that after discretization we observe unstable behavior. This phenomenon may have its roots in the non-self-adjointness of the continuous system and at this time it is unclear to us how a discretization of such a system changes its behavioral pattern.

## 9. Conclusions

In this paper we studied the dynamics of the longitudinal vibrations of a bar subjected to viscous boundary conditions at each end. The system is not self-adjoint which prevents readily obtaining the orthogonality condition. Nevertheless, it was possible to deploy Fourier's method and to derive the complete solution of the problem in terms of eigenfunction expansion. The Fourier coefficients were calculated by extending the space of the wave operator and recasting the underlying eigenvalue problem into the first order form. In the process we verified the known cases of longitudinal vibrations of a free-free, fixed-damper, fixed-fixed and fixed-free bar. We also presented how the analytical solution compares to FEM approximation and demonstrated its efficiency.

Longitudinal vibrations of bar with discrete dampers attached represent only one way to obtain non-self-adjointness in the system. Attaching dampers to any other continuous structures such as beams or plates would introduce a similar situation with greater level of complexity. To the best of our knowledge such investigations have not been carried out in the engineering literature and only some attempts to obtain approximate solution were made. However, the methodology described in this paper should in general be applicable and a similar treatment, in case of beam or plates, should yield the series form solutions. The cases of transversal vibrations of a clamped beam with a shearing damper as well as the case of transversal vibrations of a clamped beam with a torsional damper at the boundary were analyzed in [9] and [10] respectively.

## Appendix: Summation formulae

**Summation formula for** $\sum_{r=-\infty}^{\infty} \dfrac{u_{1,r}(L)u_{1,r}(x)}{\lambda_r}$

Here we determine several sums that are needed to simplify the vibratory response. We first expand a vector function $\mathbf{F}(x) = (f(x), 0)^T$ in the eigenfunctions of the operator $\mathbf{T}$ by using Eqs. (27) and (28) where the adjoint eigenvector is given by Eq. (26) to obtain

$$\begin{pmatrix} f(x) \\ 0 \end{pmatrix} = \frac{1}{L(1-h_1^2)} \sum_{r=-\infty}^{\infty} \int_0^L f(\xi)\, u_{1,r}(\xi) d\xi \begin{pmatrix} u_{1,r}(x) \\ u_{2,r}(x) \end{pmatrix}. \tag{A.1}$$

This series converges uniformly in both components, see §9.2, Theorem 4 in [7]. This provides us with the Fourier coefficients needed in Eq. (41). It follows that

$$f(x) = \frac{1}{L(1-h_1^2)} \sum_{r=-\infty}^{\infty} \int_0^L f(\xi)\, u_{1,r}(\xi) d\xi\, u_{1,r}(x) \tag{A.2}$$

$$0 = \frac{1}{L(1-h_1^2)} \sum_{r=-\infty}^{\infty} \int_0^L f(\xi)\, u_{1,r}(\xi) d\xi\, u_{2,r}(x).$$

Note that these are not Fourier series expansions in bases $u_{1,r}(x)$ and $u_{2,r}(x)$, because they are not minimal systems. Instead they are 2-fold complete [11], meaning, roughly speaking, that each function can be expanded over just "half" of the functions in the system. We also note that Parseval's theorem in this case becomes $<\mathbf{F}(x), \mathbf{F}(x)> = \int_0^L f(x)^2\, dx = \sum_{r=-\infty}^{\infty} \left( \int_0^L f(\xi)\, u_{1,r}(\xi) \right)^2 d\xi \Big/ L(1-h_1^2)$.

To obtain further simplification of the general response for Eq. (45) we set $f(x) = 1$ in (A.3) which yields

$$\begin{pmatrix} 1 \\ 0 \end{pmatrix} = \frac{1}{L(1-h_1^2)} \sum_{r=-\infty}^{\infty} \int_0^L u_{1,r}(\xi) d\xi \begin{pmatrix} u_{1,r}(x) \\ u_{2,r}(x) \end{pmatrix}. \tag{A.3}$$

From Eq. (17) we have that $u_{1,r}(x) = \dfrac{c^2 u'_{2,r}(x)}{\lambda_r}$, and in view of boundary conditions $\left(u_{1,r}(0) = 1\right)$ we have

$$\int_0^L u_{1,r}(\xi) d\xi = \frac{c^2}{\lambda_r} \int_0^L u'_{2,r}(x) dx = \frac{c^2}{\lambda_r} \left[u_{2,r}(x)\right]_0^L = -\frac{c}{\lambda_r}\left(h_2 u_{1,r}(L) + h_1\right). \tag{A.4}$$

The right-hand side of Eq. (A.4) can be written as $-\dfrac{c}{\lambda_r}\left[(h_1+h_2)u_{1,r}(L) + h_1\left(1 - u_{1,r}(L)\right)\right]$. Substituting this into Eq. (A.3) and writing out the first component of the vector equation we obtain

$$\frac{1}{(h_1+h_2)} = -\frac{c}{L(1-h_1^2)}\sum_{r=-\infty}^{\infty}\frac{u_{1,r}(L)}{\lambda_r}u_{1,r}(x) + \frac{c\,h_1}{L(h_1+h_2)(1-h_1^2)}\sum_{r=-\infty}^{\infty}\frac{u_{1,r}(L)-1}{\lambda_r}u_{1,r}(x). \tag{A.5}$$

We now want to show that the second sum on the right-hand side of Eq. (A.5) is zero. As before, we proceed with the generalized Fourier series. In view of $u_{2,r}(x) = \dfrac{u'_{1,r}(x)}{\lambda_r}$ from Eq. (17) we write

$$\begin{pmatrix}0\\1\end{pmatrix} = -\frac{c^2}{L(1-h_1^2)}\sum_{r=-\infty}^{\infty}\int_0^L u_{2,r}(\xi)d\xi\begin{pmatrix}u_{1,r}(x)\\u_{2,r}(x)\end{pmatrix} = -\frac{c^2}{L(1-h_1^2)}\sum_{r=-\infty}^{\infty}\int_0^L \frac{u'_{1,r}(x)}{\lambda_r}d\xi\begin{pmatrix}u_{1,r}(x)\\u_{2,r}(x)\end{pmatrix}. \tag{A.6}$$

Writing out the first component of the vector Eq. (A.6) we obtain

$$0 = -\frac{c^2}{L(1-h_1^2)}\sum_{r=-\infty}^{\infty}\frac{u_{1,r}(L)-u_{1,r}(0)}{\lambda_r}u_{1,r}(x). \tag{A.7}$$

Therefore Eq. (A.3) reduces to

$$\frac{1}{c(h_1+h_2)} = -\sum_{r=-\infty}^{\infty}\frac{u_{1,r}(L)u_{1,r}(x)}{\lambda_r(1-h_1^2)L}. \tag{A.8}$$

**Summation formula for** $\displaystyle\sum_{r=-\infty}^{\infty}\frac{u_{1,r}(\xi)u_{1,r}(x)}{\lambda_r}$

We now want to show that the restriction $u_{1,r}(L)$ in Eq. (A.8) can be relaxed to $u_{1,r}(\xi)$. Once again, we proceed with the generalized Fourier series by setting $f(x) = \delta(x-\xi)$ in (A.3) which yields

$$\begin{pmatrix}\delta(x-\xi)\\0\end{pmatrix} = \frac{1}{L(1-h_1^2)}\sum_{r=-\infty}^{\infty}u_{1,r}(\xi)\begin{pmatrix}u_{1,r}(x)\\u_{2,r}(x)\end{pmatrix}. \tag{A.9}$$

In view of $u_{2,r}(x) = \dfrac{u'_{1,r}(x)}{\lambda_r}$ from Eq. (17) we now write the second component of Eq. (A.9) as

$$0 = \frac{1}{L(1-h_1^2)}\sum_{r=-\infty}^{\infty}\frac{u_{1,r}(\xi)u'_{1,r}(x)}{\lambda_r}. \tag{A.10}$$

The first and second components on the right-hand side of Eq. (A.9) are not usual series, because Dirac's delta function is not a function in the ordinary sense. These series do not converge to the left-hand side of Eq. (A.9) in the classical sense, but oscillate faster and faster as more terms are added. The series are in fact distributions and therefore they converge distributionally or "weakly", see section 3.2 in [12]. Therefore, the series can be used in ordinary sense for term-by-term integration. We utilize this convergence in integrated sense and proceed by integrating Eq. (A.10) on both sides with respect to *x* to obtain

$$C(\xi) = \frac{1}{L(1-h_1^2)}\sum_{r=-\infty}^{\infty}\frac{u_{1,r}(\xi)u_{1,r}(x)}{\lambda_r}. \tag{A.11}$$

$C(\xi)$ is a constant of integration combining the constants of integration from both side of Eq. (A.11). To determine it, we set $x = L$ in Eq. (A.11) and $x = \xi$ in Eq. (A.8) after which we compare these two equations and conclude that $C(\xi) = -\dfrac{1}{c(h_1 + h_2)}$. Finally, we obtain a generalized summation formula

$$\frac{1}{c(h_1+h_2)} = -\sum_{r=-\infty}^{\infty} \frac{u_{1,r}(\xi)u_{1,r}(x)}{\lambda_r(1-h_1^2)L} \tag{A.12}$$

which can be used to simplify the response of the system.